\documentclass[12pt,letterpaper,titlepage]{amsart}
\usepackage{amsmath, amssymb, amsthm, amsfonts,amscd,xr}

\newcommand{\Z}{\mathbb{Z}}
\newcommand{\R}{\mathbb{R}}
\newcommand{\C}{\mathbb{C}}

\def\beq{\begin{equation}}
\def\eeq{\end{equation}}

\def\arr{\hbox to 20pt{\rightarrowfill}}

\def\e{\varepsilon}

\def\sl{\mathfrak {sl}}
\def\so{\mathfrak {so}}

\def\H{\mathcal H}

\newenvironment{res}
               {\begin{equation}
\begin{minipage}{0.85\textwidth}}
               { \end{minipage}\end{equation} }
\def\ber{\begin{res} }
\def\eer{\end{res}}

\numberwithin{equation}{section}
\newtheorem{thm}{Theorem}[section]

\newtheorem{lem}[thm]{Lemma}

\newtheorem{cor}[thm]{Corollary}

\newtheorem{prop}[thm]{Proposition}

\newtheorem{rem}[thm]{Remark}
\newtheorem{problem}[thm]{Problem}

\makeatletter
\def\section{\@startsection {section}{1}{\z@}{3.5ex plus 1ex minus
    .2ex}{2.3ex plus .2ex}{\large\bf}}
    \def\subsection{\@startsection{subsection}{2}{\z@}{3.25ex plus 1ex minus
 .2ex}{1.5ex plus .2ex}{\bf}}
\makeatother
\def\pf{{\em Proof}.\, }

\def\e{\epsilon}

\def\af{\mathfrak{a}}
\def\gf{\mathfrak{g}}
\def\qf{\mathfrak{q}}
\def\h{\mathfrak{h}}

\def\kf{\mathfrak{k}}
\def\lf{\mathfrak{l}}
\def\m{\mathfrak{m}}
\def\nf{\mathfrak{n}}

\def\pf{\mathfrak{p}}

\def\spec{\operatorname{spec}}

\def\zf{\mathfrak{z}}

\def\Ad{\operatorname{Ad}}
\def\ad{\operatorname{ad}}

\def\Nc{\mathcal{N}}

\def\Cc{\mathcal{C}}

\def\H{\mathcal{H}}

\def\J{\mathcal{J}}

\def\PP{\mathbb{P}}

\def\Sl{\operatorname{Sl}}

\def\W{\mathcal{W}}

\renewcommand{\Im}{\mbox{\rm Im}\,}
\makeindex
\begin{document}

\title[Domains of holomorphy]
{Domains of holomorphy for irreducible unitary representations of 
simple Lie groups}

\author{Bernhard Kr\"otz}
\address{Max-Planck-Institut f\"ur Mathematik, Vivatsgasse 7,
D-53111 Bonn \ email: {\tt kroetz@mpim-bonn.mpg.de}}

\date{\today}
\maketitle

\section{Introduction}
Let us consider a unitary irreducible representation $(\pi, \H)$
of a simple, non-compact and connected algebraic Lie group $G$.
Let us denote by $K$  a maximal compact subgroup of $G$. 
According to Harish-Chandra, the Lie algebra submodule $\H_K$
of $K$-finite vectors of $\pi$ consists of analytic 
vectors for the representation, i.e. for all $v\in \H_K$ the orbit map
$$f_v: G\to\H, \ \ g\mapsto \pi(g)v$$
is real analytic. For these functions $f_v$
we determine, and in full generality,  their natural domain
of definition as holomorphic functions (see Theorem \ref{mt} below):

\begin{thm} Let $(\pi, \H)$ be a unitary irreducible representation 
of $G$. Let $v\in \H$ be a non-zero $K$-finite vector and 
$f_v$ be the corresponding orbit map. Then there exists
a unique maximal $G\times K_\C$-invariant domain 
$D_\pi\subseteq G_\C$, independent of $v$, 
to which $f_v$ extends holomorphically. 
Explicitly: 

\begin{enumerate}
\item $D_\pi=G_\C$ if $\pi$ is the trivial representation. 
\item $D_\pi=\Xi^+ K_\C $ if $G$ is Hermitian and $\pi$ is a 
non-trivial highest weight representation. 
\item $D_\pi=\Xi^- K_\C$ if $G$ is Hermitian and $\pi$ is a 
non-trivial lowest weight representation. 
\item $D_\pi=\Xi K_\C$ in all other cases. 
\end{enumerate}
\end{thm}

In the theorem above $\Xi, \Xi^+,\Xi^-$ 
are certain $G$-domains in $X_\C=G_\C/K_\C$ over $X=G/K$
with proper $G$-action. 
These domains are studied in this paper because 
of their relevance for the theorem above (see \cite{KO}). 
Let us mention that $\Xi$ is the familiar crown domain and that 
the inclusion $\Xi K_\C\subset D_\pi$ traces back to our 
joint work with Robert Stanton (\cite{KSI}, \cite{KSII}). 

\bigskip\noindent  {\it Acknowledgment:} I am happy to point out that this 
paper is related to joint work with  Eric M. Opdam \cite{KO}. 
Also I would like to thank Joseph Bernstein who, over the years, 
helped me with his comments to understand the material much better. 
\par Finally I appreciate the work of a very good referee who
made many useful remarks on style and organization 
of the paper.

\section{Notation}

Throughout this paper $G$ shall denote a connected simple 
non-compact Lie group.
We denote by $G_\C$ the universal complexification 
of $G$ and suppose: 

\begin{itemize} \item $G\subseteq G_\C$; 
\item $G_\C$ is simply connected. 
\end{itemize}

We fix a maximal compact subgroup $K<G$ and form 
$$X=G/K\, , $$
the associated Riemannian symmetric space of the non-compact type. 
The universal complexification $K_\C$ of $K$ will be realized
as a subgroup of $G_\C$. We set 

$$X_\C=G_\C/K_\C$$
and call $X_\C$ the {\it affine complexification} of $X$. 
Note that 

$$X\hookrightarrow  X_\C, \ \ gK\mapsto gK_\C$$ 
defines a $G$-equivariant embedding which realizes
$X$ as a totally real form of the Stein symmetric space 
$X_\C$. We write $x_0=K_\C\in X_\C$ for the standard base point 
in $X_\C$.

\par However, the natural complexification of $X$ is not 
$X_\C$, but the {\it crown domain} $\Xi\subsetneq X_\C$ whose 
definition we recall now. We shall provide the standard definition of  
$\Xi$, see \cite{AG}. 
\par Lie algebras of subgroups $L<G$ will be denoted by the corresponding 
lower case German letter, i.e. ${\mathfrak l} < \gf$; 
complexifications of Lie algebras are marked 
with a $\C$-subscript, i.e. $\lf_\C$ is the complexification of $\lf$. 
\par Let us denote by $\pf$ the orthogonal complement to $\kf$ in 
$\gf$ with respect to the Cartan-Killing form. We set 

$$\hat\Omega=\{ Y\in\pf\mid \spec (\ad Y)\subset (-\pi/2,
\pi/2)\}\, .$$
Then 
$$\Xi=G\exp(i\hat \Omega)\cdot x_0\subset X_\C$$
is a $G$-invariant neighborhood of $X$ in $X_\C$, commonly 
referred to as {\it crown domain}.  
Sometimes it is useful to have an alternative, although 
less invariant picture of the crown domain: if $\af\subset\pf$ 
is a maximal abelian subspace and $\Omega:=\hat\Omega\cap \pf$, then 
\begin{equation}\label{e}\Xi=G\exp(i\Omega)\cdot x_0\, .
\end{equation}
The set $\Omega$ is nicely described through the restricted 
root system $\Sigma=\Sigma(\gf,\af)$:  
$$\Omega=\{ Y\in \af\mid \alpha(Y)< \pi/2\ \forall \alpha\in\Sigma \}\, .$$ 
If $\W$ is the Weyl group of $\Sigma$, then we note that $\Omega$
is $\W$-invariant. 

\par Sometimes we will employ the root space decomposition
$\gf=\af \oplus \m \oplus\bigoplus_{\alpha\in\Sigma}\gf^\alpha$
with $\m=\zf_\kf(\af)$ as usual. We choose a positive 
system $\Sigma^+\subset \Sigma$ and form the nilpotent subalgebra 
$\nf=\bigoplus_{\alpha\in\Sigma^+} \gf^\alpha$. 

\subsection{The example of $G=\Sl(2,\R)$} For illustration and later use  
we will exemplify the above notions at the basic 
case of $G=\Sl(2,\R)$. 

\par We let $K=\mathrm{SO}(2,\R)$ be our choice for the maximal compact
subgroup  and identify $X=G/K$  with the upper half plane 
$D^+:=\{z\in \C\mid \Im z>0\}$. 
We recall that 
$$X_\C=\PP^1(\C)\times \PP^1(\C)\setminus \mathrm{diag} [\PP^1(\C)]$$
with $G_\C$ acting diagonally by fractional 
linear transformations. The $G$-embedding of $X=D^+$ into $X_\C$ is 
given by 
$$z\mapsto (z, \overline z)\in X_\C\, .$$ 
If $D^-$ denotes the lower half plane, then the  crown domain is 
given by
$$\Xi=D^+\times D^- \subseteq X_\C\, .$$ 

\par In addition we record two $G$-domains in 
$X_\C$ which sit above $\Xi$, namely: 
\begin{align}\label{plus}\Xi^+&=D^+\times \PP^1(\C)\setminus 
\mathrm{diag} [\PP^1(\C)]\, ,\\ 
\label{min}\Xi^-&=\PP^1(\C)\times D^-\setminus \mathrm{diag} [\PP^1(\C)]\, .
\end{align}
Observe that $\Xi=\Xi^+\cap\Xi^-$.

\section{Remarks on $G$-invariant domains in $X_\C$ with proper action}

\par One defines elliptic elements in $X_\C$ by 
$$X_{\C, {\rm ell}}=G\exp(i\pf)\cdot x_0=G\exp(i\af)\cdot x_0\, .$$
The main result of \cite{AG} was to show that 
$\Xi$ is a maximal domain in $X_{\C,{\rm ell}}$ with 
$G$-action proper. In particular,  $G$ acts properly on $\Xi$.  
\par It was found in \cite{KO} that $\Xi$ in general is not a maximal domain in $X_\C$
for proper $G$-action: the domains $\Xi^+$ and $\Xi^-$ from 
(\ref{plus})-(\ref{min}) yield 
counterexamples. 
To know all maximal domains is important 
for the theory of representations \cite{KO}, Sect. 4. 

\par That $\Xi$ in general is not maximal 
for proper action is related to the unipotent model 
for the crown which was described in \cite{KO}. 
To be more precise, we showed 
that  there exists a 
domain $\hat\Lambda\subseteq \nf$ containing $0$ such that  
\begin{equation}\label{u} \Xi=G\exp(i\hat\Lambda)\cdot x_0\, .\end{equation}

Now there is a big difference between the unipotent 
parametrization (\ref{u}) and the elliptic parametrization 
(\ref{e}): If we enlarge $\Omega$ the result is no longer 
open; in particular,  $X_{\C, {\rm ell}}$ is not a domain. 
On the other hand, if we enlarge  the open set $\hat \Lambda$ the resulting set 
is still open; in particular $X_{\C, {\rm u}}:= G\exp (i\nf)\cdot x_0$ 
is a domain. 
Thus, if there were a bigger domain 
than $\Xi$ with proper action, then it is likely 
by enlargement of $\hat\Lambda$. 
\par We need some facts on the boundary of $\Xi$.

\subsection{Boundary of $\Xi$}
Let us denote by $\partial \Xi$ the topological 
boundary of $\Xi$ in $X_\C$. 
One shows 
that 
$$\partial_{\rm ell} \Xi:=G\exp(i\partial \Omega)\cdot x_0\subseteq \partial \Xi$$
(cf.\ \cite{KSII}) and calls $\partial_{\rm ell}\Xi$ the {\it elliptic part} of 
$\partial \Xi$. 
We define the {\it unipotent part} $\partial_{\rm u} \Xi $ of $\partial\Xi$ to be the 
complement to the elliptic part: 

$$\partial_{\rm u}\Xi =\partial \Xi \setminus \partial_{\rm ell}\Xi\, .$$
The relevance of $\partial_{\rm u}\Xi$ is as follows. 
Let $X\subset  D\subseteq X_\C$ denote 
a $G$-domain with proper $G$-action. Then 
$D\cap \partial_{\rm ell}\Xi=\emptyset$ 
by the above cited result of \cite{AG}. 
Thus if $D\not \subset \Xi$, then  one has 
$$D\cap \partial_{\rm u}\Xi\neq \emptyset\, .$$
\par Let us describe $\partial_{\rm u}\Xi $ in more detail. 
For $Y\in \af$ we 
define a reductive subalgebra of $\gf_\C$ by 
$$\gf_\C[Y]=\{ Z\in \gf_\C\mid e^{-2i \ad(Y)} \circ \sigma (Z)=Z\}$$
with $\sigma$ the Cartan involution on $\gf_\C$ which fixes 
$\kf + i\pf$. 
Then there is a partial result on $\partial_{\rm u}\Xi$, 
for instance stated in \cite{FH}: 

\begin{align}\label{bu}\partial_{\rm u} \Xi\subseteq\{ & G\exp(e)\exp(iY)\cdot x_0\mid 
Y\in \partial\Omega,\\
& \ 0\neq e \in \gf_\C[Y]\cap i\gf\ \hbox{nilpotent}\}\, .\end{align}

If $Y$ is such that only one root, say $\alpha$, attains the value 
$\pi/2$, then we call $Y$ and as well the elements in the boundary orbit 
$G\exp(e)\exp(iY)\cdot x_0$ {\it regular}. Accordingly we define the {\it regular unipotent}
boundary $\partial_{\rm u, reg}\Xi =\{ z\in \partial_{\rm u}\Xi\mid z\ \hbox{regular}\}$. 
Note that $\gf_\C[Y]$ is of especially simple form for regular 
$Y$, namely 

$$\gf_\C[Y]= i \af \oplus\mathfrak{m} \oplus \gf[\alpha]^{-\theta} 
\oplus i\gf[\alpha]^\theta\, $$
where $\gf[\alpha]=\gf^\alpha \oplus \gf^{-\alpha}$. Hence, in the regular situation,  
one can choose $e$ above to be in $i\gf[\alpha]^{\theta}+i\af$. 
We summarize our discussion: 

\begin{prop} \label{prop=cD}Let $X\subset D\subseteq X_\C$ be a 
$G$-invariant domain with proper $G$-action which is not contained 
in $\Xi$. Then $D\cap \partial_{\rm u, reg}\Xi\neq \emptyset$.
 More precisely, there exists $Y\in \partial \Omega$ regular (with $\alpha\in \Sigma$ 
the unique root attaining $\pi/2$ on $Y$) and 
a non-zero nilpotent element $e\in i \gf[\alpha]^\theta +i \af$
such that 
$$\exp(e)\exp(iY)\cdot x_0\in \partial_{\rm u, reg}\Xi\cap D\, .$$
\end{prop}

\section{Maximal domains for proper action}

The aim of this section is to classify all maximal 
$G$-domains in $X_\C$ which contain $X$ and maintain proper 
action. The answer will depend whether $G$ is of Hermitian 
type or not. 

\subsection{Non-Hermitian groups.}

The objective is to prove the following theorem:

\begin{thm}\label{tnh} Suppose that $G$ is not of Hermitian type. If 
$X\subset D\subset X_\C$ is a $G$-invariant domain with 
proper $G$-action, then $D\subset\Xi$. 
\end{thm}

Before we can give the proof of the theorem some preparation 
is needed. The proof relies partly on a structural fact characterizing 
non-Hermitian groups (see Lemma \ref{nh} below)
and on a precise knowledge of the basic case 
of $G=\Sl(2,\R)$.

\par Let us begin with the relevant facts for $G=\Sl(2,\R)$.  
With $E=\begin{pmatrix} 0 & 1\\ 0 & 0\end{pmatrix}$ and 
$T=\begin{pmatrix} 1& 0\\ 0 & -1\end{pmatrix}$ our choices for
$\af$ and $\nf$ are 

$$\af=\R \cdot T \quad \hbox{and}\quad \nf=\R\cdot E\, .$$
Note that $\Omega=(-\pi/4,\pi/4)T$.

\par The a slight modification of results in  \cite{KO}, Sect. 3 and 4
yield: 

\begin{lem}\label{pa1} Let $G=\Sl(2,\R)$ and $\J\subset\R$ be an 
open subset. Then 
$$\Xi_\J:=G\exp(iJ\cdot E)\cdot x_0$$
is a $G$-invariant open subset of $X_\C$ and the following holds: 
\begin{enumerate}
\item $G$ does not act properly if $\{-1,1\}\subset\J$. 
\item $\Xi=\Xi_{(-1,1)}$. 
\item $\Xi^+=\Xi_{(-1,\infty)}$. 
\item $\Xi^-=\Xi_{(-\infty, 1)}$.  
\end{enumerate}
\end{lem}

We also need that $\partial \Xi$ is a fiber bundle over 
the affine symmetric space $G/H$ where $H=\mathrm{SO}_e(1,1)$.
Notice that $H$ is the stabilizer of the 
boundary point 

$$z_H:=\exp(-i\pi T/4)\cdot x_0=(1,-1)\in \partial_{\rm ell}\Xi\, .$$

Write $\tau$ for the involution on $G$, resp. $\gf$, fixing 
$H$, resp. $\h$,  and denote by 
$\gf=\h+\qf$ the corresponding eigenspace decomposition. 
The $\h$-module $\qf$ breaks into two eigenspaces 
$\qf=\qf^+ \oplus \qf^-$ with 

$$\qf^\pm =\R \cdot e^\pm \quad \hbox{where}\quad 
e^\pm =\begin{pmatrix} 1 & \mp 1 \\ \pm 1 & - 1\end{pmatrix}\, .$$
Finally write
$$\Cc=\R_{\geq 0}\cdot e^+ \cup \R_{\geq 0}\cdot e^-$$
and $\Cc^\times =\Cc\setminus\{0\}$. 
Note that both $\Cc$ and $\Cc^\times$ are $H$-stable.  
We cite \cite{KO}, Th. 3.1: 

\begin{lem}\label{pa3} Let $G=\Sl(2,\R)$. Then the map 
$$G\times_H \Cc\to \partial \Xi, \ \ [g,e]\mapsto g\exp(ie)\cdot z_H$$
is a $G$-equivariant homeomorphism.  
Moreover, 
\begin{enumerate}
\item $\partial_{\rm ell}\Xi=G\cdot z_H\simeq G/H$, 
\item $\partial_{\rm u}\Xi=G\exp(i\Cc^\times)\cdot 
z_H\simeq G\times_H \Cc^\times$, 
\item $\partial_{\rm u}\Xi =G\exp(i E)\cdot x_0\amalg G\exp(-i E)\cdot x_0 $. 
\end{enumerate}
\end{lem}

As a last piece of information we need a structural fact which is  
only valid for non-Hermitian groups. 

\begin{lem}\label{nh} Suppose that $G$ is not of Hermitian type. 
Then for all $\alpha\in \Sigma$ and 
$E\in\gf^\alpha$ there exists an $m\in M=Z_K(\af)$ such that 
$$\Ad(m) E = -E\, .$$
\end{lem}

\begin{proof} Let us remark first that we may assume that $G$ is of 
adjoint type. If $G$ is complex, then the assertion is clear 
as $T:=\exp(i\af)\subset M$ provides us with the elements 
we are looking for. More generally for $\dim \gf^\alpha>1$ 
one knows (Kostant) that $M_0=\exp(\m)$ acts 
transitively on the unit sphere in $\gf^\alpha$ (cf. \  \cite{Kos}). 
\par In the sequel we use the terminology and tables of the classification 
of real simple Lie algebras as found in the monograph \cite{K}, App. C. 
As $G$ is not Hermitian, Kostant's result leaves us with the 
following cases for $\gf$: $\sl(n,\R)$ for $n\geq 3$, 
$\so(p,q)$ for $0,2\neq p,q$ and $p+q>2$, 
$E\ I$, $E\ II$, $E\ V$, $E\ VI$, $E\ VIII$, $E\ IX$, $F\ I$ and $G$.  
\par Now we make the following observation. The lemma is true 
for $G=\Sl(3,\R)$ as a simple matrix computation shows. 
Suppose that $\alpha$ is such that it can be put 
into an $A_2$-subsystem of $\Sigma$. As $\dim \gf^\alpha$ is 
one-dimensional (by our reduction) this means that we can put 
$E\in \gf^\alpha$ in a subalgebra isomorphic to 
$\sl(3,\R)$. Now it is important to recall the nature of the 
component group of $M$, see \cite{K}, Th. 7.55. It follows that 
the $M$-group of $\Sl(3,\R)$ (isomorphic to $(\Z/2\Z)^2$) 
embeds into the $M$-group of $G$. 
\par The $A_2$-reduction described above deletes most of the 
cases in our list. We remain with the orthogonal cases 
$\so(p,q)$ for $0, 2 \neq p,q$ and $p\neq q$. 
A simple matrix computation, which we leave to the reader, 
finishes the proof. 
\end{proof}

\begin{proof} (of Theorem \ref{tnh}) Suppose that $G$ is not of Hermitian 
type. Let $X\subset D\subset\Xi$ be a $G$-invariant 
domain with proper $G$-action which is not contained 
in $\Xi$. We shall show that $D$ does not exist. 

\par According to Proposition \ref{prop=cD} we find a regular 
$Y\in\partial\Omega $ and a non-zero nilpotent $e\in \gf_\C[Y]\cap i\gf$ 
such that 
$$\exp(e)\exp(iY)\cdot x_0\in \partial_{\rm u, reg} \Xi\cap D\, .$$ 
Let $\alpha\in\Sigma$ be the root corresponding to $Y$. 
Write $Y=Y^\alpha +Y'$ with $Y^\alpha,Y'\in\af$ such that $\alpha(Y')=0$. 
It is known that $Y^\alpha \in\partial \Omega$ and $Y'\in\Omega$.  
Hence we may use 
$\sl(2)$-reduction which in conjunction with Lemma \ref{pa3}
implies the existence of $E^\alpha\in\gf^\alpha$ 
such that: 
\begin{itemize} 
\item $\{E^\alpha, \theta(E^\alpha), [E^\alpha, \theta(E^\alpha)]\}$ is an 
$\sl(2)$-triple, 
\item $\exp(iE^\alpha)\exp(iY')\cdot x_0 \in \partial_{\rm u, reg}\Xi\cap D$, 
\end{itemize}
Now, as $G$ is not of Hermitian type, Lemma \cite{nh} implies 
that there 
exists an element $m\in M$ such that 
$\Ad(m) E^\alpha =-E^\alpha$. 
Hence 
 $$\exp(-iE^\alpha)\exp(iY')\cdot x_0\in \partial_{\rm u, reg}\Xi$$
as well. But this contradicts Lemma \ref{pa1}(i). 
\end{proof}

\subsection{Hermitian groups}

Let  now $G$ be of Hermitian type and $G\subseteq P^- K_\C P^+$ be a Harish-Chandra
decomposition of $G$ in $G_\C$. We define flag varieties 
$$F^+=G_\C/K_\C P^+\quad\hbox{and}\quad F^-=G_\C/ K_\C P^-$$
and inside of them we declare the flag domains
$$D^+ =GK_\C P^+ /K_\C P^+\quad\hbox{and}\quad D^-=GK_\C P^-/ K_\C P^-\, .$$
Then 
\begin{equation} \label{em1} 
X_\C \hookrightarrow F^+\times F^-,\ \  gK_\C \mapsto 
(gK_\C P^+, gK_\C P^-) \end{equation}
identifies $X_\C$ as a Zariski open affine piece of 
$F^+\times F^-$. In more detail: As $G$ is of Hermitian type, 
there exist $w_0\in N_{G_\C}(K_\C)$ 
such that $w_0P^\pm w_0^{-1} = P^{\mp}$. In turn, this element induces 
a $G_\C$-equivariant biholomorphic map: 

$$\phi: F^+ \to F^-, \ \ gK_\C P^+\mapsto gw_0 K_\C P^-\, .$$
With that the embedding (\ref{em1}) gives the following 
identification of $X_\C$:  
\begin{equation} \label{xc} 
X_\C=\{ (z,w)\in F^+\times F^-\mid \phi(z)\intercal w\}\, , \end{equation}
where $\intercal$ stands for the transversality notion in the 
flag variety $F^-$. 
We recall what it means to be transversal. First note that the notion 
is $G_\C$-invariant, i.e. for $z,w\in F^-$ and $g\in G_\C$ one has 
$z\intercal w$ if and only if  $gz\intercal gw$. Now for the base point 
$z^-= K_\C P^-\in F^-$ one has $z^-\intercal w$ 
if and only if $w\in P^-w_0 z^-$.  

\par  We keep  the realization of $X_\C$ 
in $F^+\times F^-$ (cf. (\ref{em1})
in mind and recall the description of $\Xi$: 
$$\Xi=D^+\times D^-\, $$
(see \cite{KSII}). 
\par For subsets $X^\pm \subset F^\pm$ we write 
$X^+\times_\intercal X^-$ for those elements $(x^+, x^-)\in X^+\times X^-$
which are transversal, i.e. $\phi(x^+)\intercal x^-$. 
With this terminology in mind we finally define 
$$\Xi^+= D^+ \times_\intercal F^- ,$$
$$\Xi^-= F^+\times_\intercal D^-\, .$$

\subsubsection{Basic structure theory of $\Xi^+$ and $\Xi^-$}

It is obvious that both $\Xi^+$ and $\Xi^-$ are open and $G$-invariant. 
However, as was pointed out by the referee,  
it is a priori not clear that they are connected. 
In order to see this let $p_+: \Xi^+ \to D^+$ 
be the projection onto the first factor. Likewise 
we define $p_-:\Xi^-\to D^-$.

\begin{prop}\label{p17} Let $\e\in \{-, +\}$. The map $p_\e: \Xi^\e\to D^\e$ 
induces the 
structure of a holomorphic fiber bundle with fiber 
isomorphic to $P^\e$. 

\end{prop}

\begin{proof} We confine ourselves with the case $\e=+$. 
\par As $p_+$ is $G$-equivariant and $D^+$ is $G$-homogeneous, 
it is sufficient to determine the fiber $p_+^{-1}(z^+)$. 
Recall that  $z^+=K_\C P^+\in F^+$ is the base point. 
Now 
$$p_+^{-1}(z^+)=\{ (z^+, w)\in F^+\times F^-\mid \phi(z^+)
\intercal w\}\, .$$ 
Observe that $\phi(z^+)=w_0 z^-$ and that 
$w_0z^- \intercal w $ is equivalent to $z^-\intercal w_0^{-1}w$. 
By the definition of transversality this means that 
$w_0^{-1} w \in P^- w_0 z^-$ or $w\in w_0 P^- w_0 z^-$. 
It is no loss of generality to assume that $w_0=w_0^{-1}$. 
So we arrive at $w\in P^+ z^-$ and this concludes 
the proof of the proposition.
\end{proof}

\begin{cor} Both $\Xi^+$ and $\Xi^-$ are contractible. 
\end{cor}

It was observed by the the referee that Proposition \ref{p17} 
allows the following interesting reformulation.  

\begin{cor} \label{c1}The  map 
$$G\times_K P^+\to \Xi^+, \ \ [g, p]\mapsto 
(gz^+,g pz^-)$$
is a $G$-equivariant diffeomorphism. In particular 
$\Xi^+$ is $G$-biholomorphic to $T^{0,1} D^+$, the 
antiholomorphic tangent bundle of $D^+$. 
Likewise, $\Xi^-$ is $G$-biholomorphic to $T^{0,1} D^-$. 
\end{cor} 
  
Corollary \ref{c1} combined with the Harish-Chandra decomposition 
implies that 
$\Xi^\e\simeq D^\e \times P^\e$ as complex manifolds. In particular 
$\Xi^\e$ is Stein. 

\par The fact that $K_\C$ normalizes $P^\e$ allows us to speak 
of $G\times P^\e$-invariant domains in $X_\C$. 
It follows from (\ref{em1}) and Corollary \ref {c1} that $\Xi^\e$ is 
$G\times P^\e$-invariant.

\begin{prop}\label{p22} Let $\e\in\{-, +\}$. 
The real group $G$ acts properly on $\Xi^\e$. Moreover $\Xi^\e$ 
is a maximal  $G\times P^\e$-invariant domain in $X_\C$ for proper 
$G$-action.
\end{prop}

\begin{proof} As the $G$-action is proper on $D^\e$, it follows 
that $G$ acts properly on $\Xi^\e$. 
In the sequel we deal with $\e=+$ only. 
It remains to show that $\Xi^+$ is a maximal $G\times P^+$-invariant 
domain in $X_\C$ for 
proper $G$-action.
We argue by contradiction and suppose that 
$D\supsetneq\Xi^+$ is a $G\times P^+$-domain in $X_\C$ with proper 
$G$-action.
Then $D=(D_0\times F^-)\cap X_\C$ with $D_0\supsetneq D^+$ a $G$-domain 
with proper action. Now recall the following 
facts: 
\begin{itemize}
\item There are only finitely many 
$G$-orbits in $F^+$. 
\item There are precisely two orbits
with proper $G$-action: $D^+$ and $\phi^{-1}(D^-)$. 
\end{itemize}
The assertion follows. 
\end{proof}

\begin{rem} Suppose that $G$ is of Hermitian type. 
Then it can be shown that if $X\subseteq D\subseteq X_\C$ is a 
$G$-invariant domain with proper $G$-action, then 
$D\subseteq \Xi^+$ or $D\subseteq\Xi^-$. 
\par As we will not need this fact, we refrain from 
a proof. \end{rem}

If $D\subseteq X_\C$ is a subset, then we write $DK_\C$ 
for its preimage in $G_\C$ under the canonical 
projection $G_\C\to X_\C$. 

\begin{prop}\label{pi} The following assertions hold: 
\begin{enumerate}
\item $\Xi^+K_\C =GK_\C P^+$,  
\item $\Xi^-K_\C= G K_\C P^-$. 
\end{enumerate}
\end{prop}

\begin{proof} It suffices to prove (i).  Recall the embedding (\ref{em1}), 
and the definition of transversality condition. We deduce that 
$P^+\subset \Xi^+ K_\C$. As $\Xi^+K_\C$ is $G\times K_\C$-invariant, 
it follows that $GP^+K_\C=G K_\C P^+\subset \Xi^+ K_\C$. 
\par Conversely, Corollary \ref{c1} implies that $G P^+$ maps onto $\Xi^+$
and thus $\Xi^+ \subset GP^+ K_\C$. 
\end{proof}

We conclude this subsection with some easy facts on the structure of 
$\Xi^+$ and $\Xi^-$ which will be used later on.

\subsubsection{Unipotent model for $\Xi^+$ and $\Xi^-$}

We begin with the unipotent parameterization of 
$\Xi^+$ and $\Xi^-$. Some terminology is needed. 

\par According to C. Moore, $\Sigma$ is of type $C_n$ or $BC_n$. 
Hence  we find 
a subset $\{\gamma_1, \ldots, \gamma_n\}$ of long 
strongly orthogonal restricted roots.
We fix $E_j\in \gf^{\gamma_j}$ such that $\{E_j, \theta(E_j), [E_j, \theta E_j]\}$
becomes an $\sl(2)$-triple. Set $T_j:=1/2 [E_j,\theta E_j]$ and note that 
$$\Omega=\bigoplus_{j=1}^n (-\pi/2,\pi/2) T_j\, .$$ 
We set $V=\bigoplus_{j=1}^n \R\cdot E_j$ and take a
cube inside $V$ by 
$$\Lambda=\bigoplus_{j=1}^n (-1,1)E_j\, . $$
In \cite{KO}, Sect. 8,  we have shown that 
$$\Xi=G\exp(i\Lambda)\cdot x_0\, .$$
In this parametrization of $\Xi$ the unipotent boundary 
piece has a simple description: 

\begin{equation}\label{ub}\partial_{\rm u}\Xi=G\exp(i\partial \Lambda)\cdot x_0\, .\end{equation}

The strategy now is to enlarge $\Xi$ by enlarging $\Lambda$
while maintaining that the object stays a domain on which $G$ acts properly. 
But now we have to be a little bit careful 
with our choice of $E_j$. Replacing $E_j$ by $-E_j$ has no effect for the 
matters cited above, but for the sequel. Our choice is such that 
$\gamma_1, \ldots, \gamma_n$ are positive roots (this determines 
the non-compact roots in $\Sigma^+$ uniquely). 
We set 
$$\Lambda^+=\bigoplus_{j=1}^n (-1,\infty) E_j\quad \hbox{and}\quad 
\Lambda^-=\bigoplus_{j=1}^n (-\infty,1) E_j\, .$$

Then, a direct generalization of  Lemma \ref{pa1}(iii),(iv) yields: 

\begin{prop}\label{p11}The following assertions hold: 
\begin{enumerate}
\item $\Xi^+=G\exp(i\Lambda^+)\cdot x_0$, 
\item $\Xi^-=G\exp(i\Lambda^-)\cdot x_0$.
\end{enumerate}
\end{prop}

\begin{rem}  If we define subcones of 
the nilcone  $\Nc\subseteq \gf$ by 
$$\Nc^+=\Ad(K)\left[\bigoplus_{j=1}^n [0,\infty) E_j\right] 
\quad\hbox{and}\quad \Nc^-=- \Nc^+\, ,$$ 
then one can show that the maps
$$G\times_K\Nc^\pm\to\Xi^\pm, \ \ [g,Y]\mapsto g\exp(iY)\cdot x_0$$
are homeomorphic. 
\end{rem}

\section{Representation theory}

Let $(\pi,\H)$ be a unitary representation of $G$ and 
$\H_K$ the underlying Harish-Chandra module of $K$-finite 
vectors. Notice that $\H_K$ is naturally a module for 
$K_\C$. 
\par We say that $(\pi,\H)$ is a highest, resp. lowest, weight 
representation if $G$ is of Hermitian type and 
$\pf^+=\mathrm {Lie} (P^+)$, resp. $\pf^-$, acts on $\H_K$ in a finite   
manner. 

\par We turn to the main result of this paper. 

\begin{thm}\label{mt} Let $(\pi, \H)$ be a unitary irreducible representation 
of $G$. Let $v\in \H$ be a non-zero $K$-finite vector and 
$$f_v: G\to \H, \ \ g\mapsto \pi(g)v$$
the corresponding orbit map. Then there exists
a unique maximal $G\times K_\C$-invariant domain 
$D_\pi\subseteq G_\C$, independent of $v$, 
to which $f_v$ extends holomorphically. 
Explicitly: 

\begin{enumerate}
\item $D_\pi=G_\C$ if $\pi$ is the trivial representation. 
\item $D_\pi=\Xi^+ K_\C $ if $G$ is Hermitian and $\pi$ is a 
non-trivial highest weight representation. 
\item $D_\pi=\Xi^- K_\C$ if $G$ is Hermitian and $\pi$ is a 
non-trivial lowest weight representation. 
\item $D_\pi=\Xi K_\C$ in all other cases. 
\end{enumerate}
\end{thm}

\begin{proof} If $\pi$ is trivial, then the assertion is clear. 
So let us assume that $\pi$ is non-trivial in the sequel. 
Fix a nonzero $K$-finite vector $v$ and consider 
the orbit map $f_v: G\to \H$. We recall the following two facts: 

\begin{itemize}
\item $f_v$ extends to a holomorphic $G$-equivariant map
$f_v: \Xi K_\C \to \H$  (see \cite{KSII}, Th. 1.1). 
\item If $D_v\subseteq G_\C $ is a $G\times K_\C$-invariant domain to which 
$f_v$ extends holomorphically, then $G$ acts properly on $D_v/K_\C$
(see \cite{KO}, Th. 4.3)
\end{itemize}

\par We begin with the case where $G$ is not of Hermitian type. 
Here the assertion follows from the bulleted items above in conjunction 
with Theorem \ref{tnh}. 

\par So we may assume for the remainder that $G$ is of Hermitian type. 
If $\pi$ is a highest weight representation, then it is clear 
that $f_v$ extends to a holomorphic map 
$GK_\C P^+\to \H$. Thus, in this case $\Xi^+K_\C =GK_\C P^+$
(cf. Proposition \ref{pi} ) is a maximal domain 
of definition for $f_v$ by Proposition \ref{p22} and the 
second bulleted item from above. 
Likewise, if $(\pi,\H)$ is a lowest weight representation, 
then $\Xi^-K_\C$ is a maximal domain of definition 
of  $f_v$. As both $\Xi^+$ and $\Xi^-$ are simply connected with 
sufficiently regular boundary, it follows that these maximal domains are in fact unique.

\par It remains to show: 

\begin{itemize}
\item If $f_v$ extends holomorphically on a domain $D\supset\Xi$ such that $D\cap [\Xi^+\setminus \Xi]\neq \emptyset$, then 
$(\pi,\H)$ is a highest weight representation. 
\item If $f_v$ extends holomorphically p 
on a domain $D\supset\Xi$ such that $D\cap [\Xi^-\setminus \Xi]\neq \emptyset$, then 
$(\pi,\H)$ is a lowest weight representation. 
\end{itemize}
 
It is sufficient to deal with the first case. 
So suppose that $f_v$ extends to a bigger domain $D$ such that 
$D\cap [\Xi^+\setminus \Xi]\neq \emptyset$. 
Taking derivatives and applying the fact that
$d\pi({\mathcal U}(\gf_\C))v=\H_K$, we see that 
$f_u$ extends to $D$ for all $u\in \H_K$.   
By Proposition \ref{prop=cD}, (\ref{ub}) and our assumption we 
find $1\leq j\leq n$ be such that $\exp(iE_j)\exp(iY)\cdot x_0\in D$
for some $Y\in\Omega$ with $\gamma_j(Y)=0$. Let $G_j<G$ be the analytic subgroup 
corresponding 
to the $\sl(2)$-triple $\{ E_j, \theta(E_j), [E_j, \theta(E_j)]\}$. 
Basic representation theory of type I-groups in conjunction 
with \cite{KO}, Th. 4.7, yields that 
$\pi\big|_{G_j}$ breaks into a direct sum of highest weight representations. 
Applying $N_K(\af)$ (which in particular permutes the $G_k$ and preserves $\H_K$) we see 
that above matters  hold for any other $G_k$ as well (note that $Y$ might change but this does not matter
as $\Omega$ is $N_K(\af)$-invariant).  
It follows that $\pi$ is a highest weight representation
and completes the proof of the theorem. 
\end{proof}

\begin{rem} The domains $\Xi$, $\Xi^+$ and $\Xi^-$ are independent 
of the choice of the connected group $G$. Accordingly, the above theorem 
holds for all simple connected non-compact Lie 
groups $G$, i.e. we can drop the assumption that $G\subseteq G_\C$
and $G_\C$ simply connected. 
\end{rem}

\begin{problem} The above theorem should hold true for all irreducible admissible 
Banach representations of $G$ under the reservation that (i) gets modified to 
: $D_\pi=G_\C$ if $\pi$ is finite dimensional. 
\end{problem}

\end{document}